\documentclass[12pt]{article}
\usepackage{latexsym}
\usepackage{epsfig}
\usepackage{verbatim}
\usepackage{amssymb}
\usepackage{amsmath}

\title{Dixmier traces and coarse multifractal analysis}
\author{K.J. Falconer and A. Samuel\\
\small{{\it Mathematical Institute, 
University of St~Andrews, North Haugh, St~Andrews,}} \\
\small{{\it Fife, KY16~9SS, Scotland }}} 

\newtheorem{theo}{Theorem}

\newtheorem{prop}[theo]{Proposition}

\newtheorem{lem}[theo]{Lemma}

\newcommand{\bbbr}{\mathbb R}

\newcommand\dimh{{\rm dim_H}}
\newcommand\dimb{{\rm dim_B}}
\newcommand\dimlb{{\rm \underline{dim}_B}}
\newcommand\dimub{{\rm \overline{dim}_B}}
\newcommand\BB{{\overline B}}
\newcommand\IB{{\overline I}_n}

\newcommand\be{\begin{equation}}
\newcommand\ee{\end{equation}}

\newcommand{\Pf}{\par\noindent{\it Proof}\,}  
\newenvironment{pf}{\Pf}{\begin{flushright} $\square$ \end{flushright}}

\begin{document}
\maketitle
\begin{abstract}
\noindent 
We show how multifractal properties of a measure supported by a fractal $F\subseteq  [0,1]$ may be expressed in terms of complementary intervals of $F$ and thus in terms of spectral triples and the Dixmier trace of certain operators. For self-similar measures this leads to a noncommutative integral over $F$ equivalent to integration with respect to an auxilary multifractal measure.

 \end{abstract}


\section{Introduction}
\setcounter{equation}{0}
\setcounter{theo}{0}

This paper has two aims: to show how coarse multifractal analysis of a measure on a subset of $[0,1]$ may be expressed in terms of complementary intervals, and then to show how this may be used to encode multifractal information using spectral triples and the Dixmier trace of operators.

A  compact `fractal' subset  $F$ of $[0, 1]$ of Lebesgue measure $0$ may be regarded as the complement in $[0, 1]$ of a family $ \{I_n\}$ of countably many disjoint open intervals $I_n$ which we order by decreasing length. Relationships between the Hausdorff and Minkowski dimensions of $F$ and the asymptotic behaviour of the lengths $ \{|I_n|\}_{n=1}^{\infty}$have been known for a long time, see \cite{BT}, and more recently \cite{MMF,LandP}. Here we show that, given a measure $\mu$ supported by $F$, coarse multifractal aspects  of $\mu$, and in particular the  `$\beta(q)$' function, are related to the sequence of measures of intervals $\{\mu( \IB^a) \}_{n=1}^{\infty}$, where $\IB^a$ is the interval obtained by enlarging each $I_n$ by a constant factor $1+a$ about its centre.   Expressing multifactal features of a  measure $\mu$ in terms of the $\mu( \IB^a)$ is a multifractal analogue of studying the fractality of a set $F$ in terms of the $|I_n|$. We show that, under a certain  `lacunarity' condition, the $\beta(q)$ function that arises in this way is stable for both positive and negative $q$ and equals the usual $\beta(q)$ function for positive $q$.

Connes \cite{C} introduced a spectral triple on a fractal  $F \subseteq \bbbr$ defined in terms of the  complementary intervals $I_n$. This was studied by Guido and  Isola \cite{GI0, GI1, GI2} who showed that the Minkowski dimension of $F$ equals to the spectral dimension of the triple and the Minkowski measure of $F$, if it exists,  is given by the Dixmier trace of an operator. In certain cases this gives rise to a noncommutative integral over $F$. We develop multifractal analogues of this construction, with the multifractal $\beta(q)$ function expressed in terms of the Dixmier trace of certain operators. In particular, for self-similar measures this leads to a noncommutative integral over $F$ that is equivalent to integration with respect to a multifractal auxillary measure.

Section 2 gives a brief summary of spectral triples and singular traces. Section 3 reviews how fractal dimensions can be expressed in this context, and Section 4 develops these ideas for multifractals.


\section{Spectral triples and singular traces}
\setcounter{equation}{0}
\setcounter{theo}{0}

A {\it spectral triple} $(A,H, D)$ consists of a $C^{*}$-algebra $A$ acting faithfully on a Hilbert space $H$ and a Dirac operator $D: H \to H$.  That is, $D$ is a self-adjoint unbounded operator with compact resolvent such that 
$\{ a \in A :  \lVert [D, \xi(a)] \rVert  <  \infty \; \}$
is dense in $A$, where $\xi$ is the $*$--representation of $A$ on $H$ and for $a \in A$ and $[D, \xi(a)]$ is the commutator operator on $H$ defined by $[D, \xi(a)]h = D \xi(a)h - \xi(a)(D(h))$, $(h \in \mathcal{H}$).

For $H$ a separable Hilbert space and $T: H \rightarrow H$ we write $\lvert T \rvert = (T^{*}T)^{\frac{1}{2}}$. The singular values $\sigma_{k}(T)$ are taken in decreasing order, that is $\sigma_{k}(T)$ is the 
$k$th largest eigenvalue of $\lvert T \rvert$.

Let $B(H)$ denote the space of bounded operators on $H$ with $K(H)$ the ideal of compact operators. Recall that a {\it  trace} $\tau$ on  $B(H)$ is a positive linear functional defined on a    (two sided) ideal $I \subseteq K(H)$ that is unitary invariant. We will construct an operator such that the infinitesimal properties of a measure $\mu$ is reflected in the asymptotic behaviour of the singular values. This asymptotic behaviour is in turn described by a singular trace of the operator.

A {\it singular trace} $\tau: I \to \bbbr$ is a trace defined on an ideal $ I \subseteq K(H)$  such that 
if $\sigma_{i}(S) \leq \sigma_{i}(T)$ for all but a finite number of $i \in \mathbb{N}$ then $\tau(S) \leq \tau(T)$.  In particular this implies that if  $S$  has finite rank then $ \tau(S) = 0$, and also that if $\lim_{i \rightarrow \infty}\sigma_{i}(S)/\sigma_{i}(T) = 1$ then $\tau(S) = \tau(T)$.

An operator $T$ on $H$ is {\it infinitesimal of order} $\alpha$ if $\sigma_{k}(T) = O(k^{-\alpha})$ as $k \rightarrow \infty$. Note that if 
$\sigma_{k}(T) =O( k^{-1})$ then $ \sum_{k = 1}^{N} \sigma_{k}(T) = O(\log N)$; this logarithmic divergence is utilised in the Dixmier trace to quantify operators that are infinitesimal of order $1$.

Let  $\omega$ be a state on $l^\infty$, that is a positive linear functional of norm $1$.
We write $\omega \text{-}\lim a_n$ for the value of $\omega$ on $(a_n) \in l^\infty$; this notation
is indicative of the fact that  $\omega$ is often a generalised (Banach) limit. 

The {\it Dixmier ideal} of the Hilbert space $H$ is 
\be
\mathcal{L}^{1,\infty}(H) \; \mathrel{:=} \; \left\{ T \in K(H) : \limsup_{N \rightarrow \infty} \; \frac{ \sum_{k = 1}^{N} \sigma_{k}(T) }{\log N} \; < \; \infty \right\} \label{dixid}
\ee
which is a two-sided ideal of $B(H)$. 
A {\it Dixmier trace}  is a positive linear functional $\tau_\omega:\mathcal{L}^{1,\infty}(H) \rightarrow \mathbb{R} $ defined by
\be
\tau_\omega(T) \; \mathrel{:=} \omega \text{-}\lim \frac{ \sum_{k = 1}^{N} \sigma_{k}(T) }{\log N},
\label{Dixtraa}
\ee
where $\omega$ is some state. Then $\tau_\omega$ is a singular trace on $\mathcal{L}^{1,\infty}(H)$ which satisfies
$$\liminf_{N \rightarrow \infty}  \frac{ \sum_{k = 1}^{N} \sigma_{k}(T) }{\log N}  \leq \tau_\omega(T) \leq \limsup_{N \rightarrow \infty}\frac{ \sum_{k = 1}^{N} \sigma_{k}(T) }{\log N}.$$

There are many possible Dixmier traces. However, if $\tau_\omega(T) $ is independent of $\omega$, that is if the limit in (\ref{Dixtraa}) exists, we say that $T$ is {\it measurable} and refer to {\it the} Dixmier trace, and define the  {\it non-commutative integral} of an operator $T \in \mathcal{L}^{1,\infty}(H)$ to be the common value
\be
-\hspace{-0.45cm}\int T \mathrel{:=} \tau_\omega(T)= \lim_{N \rightarrow \infty} \frac{\sum_{i = 1}^{N} \sigma_{i}(T)}{\log N}. \label{Dixtr}
\ee

The {\it spectral dimension} of a spectral triple $(A, H, D)$ is defined as
\be
d(A, H, D) \mathrel{:=} \inf \big\{ s \geq 0 : |D|^{-s} \in \mathcal{L}^{1,\infty}(H) \big\}= \sup \big\{ s \geq 0 : |D|^{-s} \not\in \mathcal{L}^{1,\infty}(H) \big\}.
\ee

These notions allow  us to integrate with respect to a spectral triple $(A, H, D)$. Given $a \in A$ we may define an integral of $a$ as
\[
-\hspace{-0.45cm}\int \xi(a) \lvert D \rvert^{-d(A, H, D)}
\]
where $\xi$ is the $*$-representation of $A$ over $H$, providing the non-commutative integral exists. The exponent $ -d(A, H, D)$ is the critical parameter that is required to have any chance of a non-trivial  integral.


\section{Representing fractals by spectral triples}
\setcounter{equation}{0}
\setcounter{theo}{0}

In this section we review briefly Connes's method \cite{C} of constructing a spectral triple on a fractal, see also D.~Guido and T.~Isola \cite{GI0, GI1, GI2}. 

Let $ {\mathcal{H}^{s}(E)}$ denote the $s$-dimensional Hausdorff measure of a set $E\subseteq \bbbr$ and $\dimh E$ the the Hausdorff dimension of $E$. We write $\dimlb E$ and $\dimub E$ for the lower and upper Minkowski (or box-counting) dimensions of $F$ and if $\dimlb E=\dimub E$ we term the common value the {\it Minkowski dimension} and write $\dimb E$.
Moreover, writing $E_r$ for the $r$-neighbourhood of $E$,  the $s${\it-dimensional Minkowski content} of $F$ is given by
$M_s(E) = \lim_{r \to 0} {\mathcal L}(E_r)/r^{1-s}$ if this limit exists, where ${\mathcal L}$ is Lebesgue measure, in which case we say that $E$ is $s${\it-Minkowski measurable}. In particular, if this is the case, or even if $E$ is $s${\it-almost Minkowski measurable}, that is ${\mathcal L}(E_r)\asymp r^{1-s}$ as $r \to 0$, then   $\dimb E=s$. (The symbol $\asymp$ indicates that the ratio of the two sides is bounded away from $0$ and $\infty$.)  See, for example, \cite{F1} for these definitions and their properties.

Let $F$ be a {\it fractal set}, which here we take to mean a compact subset of $[0, 1]$ of Lebesgue measure $0$ containing the points $\{0,1\}$, so that $F$  is the complement in $[0, 1]$ of a family $ \{I_n\}$ of countably many disjoint open intervals.  We order these intervals in decreasing order of their lengths $|I_n|$, so that 
$|I_1| \geq |I_2| \geq \dots > 0 $, and write  $b_{n}^{-}$ and $b_{n}^{+}$ for the left and right endpoints  of $I_n$. Of course $\sum_{n=1}^\infty |I_n| = 1$.

To obtain a spectral triple on $F$, let $A$ be the algebra of continuous complex valued functions on $F$ and let $H$ be the Hilbert space $l^{2}(\bigcup_{n \in \mathbb{N}} \{b_{n}^{-},b_{n}^{+}\})$.  Let $\xi: A \rightarrow \mbox{End}(H)$ be defined by
\[
\xi(a) \left( \begin{array}{c} x_{i} \\ y_{i} \end{array} \right)_{i \in \mathbb{N}} \mathrel{:=} \left( \begin{array}{c} a(b_{i}^{-}) \; x_{i} \\ a(b_{i}^{+}) \; y_{i} \end{array} \right)_{i \in \mathbb{N}}
\quad \mbox{ for } \quad \left( \begin{array}{c} x_{i} \\ y_{i} \end{array} \right)_{i \in \mathbb{N}}\in H, \quad a \in A
\]
and let $D : H \rightarrow H$ be
\[
D  =  \bigoplus_{n \in \mathbb{N}}  \frac{1}{|I_n|} \left( \begin{array}{cc} 0 & 1\\ 1 & 0 \end{array} \right).
\]
Then $D$ is an unbounded operator with domain a dense subset of $H$ which has a bounded compact inverse
\[
D ^{-1} \; = \; \bigoplus_{n \in \mathbb{N}} \; |I_n|\left( \begin{array}{cc} 0 & 1\\ 1 & 0 \end{array} \right).
\]
It is easy to verify the commutator condition, so  $(A,H,D)$ is a spectral triple. 
For $s>0$ the operator $|D|^{-s}$ is compact  with eigenvalues $|I_n|^s$ each of multiplicity $2$.

We now quote two results which relate the trace of  $|D|^{-s}$ to the dimensions and Hausdorff measure of $F$. These results are essentially given in \cite{GI1}.

\begin{theo}\label{tracemm}
Let $F$ be a fractal subset of $[0, 1]$ as above.  Then 
 $d(A, H, D) = \dimub F$. If ${\mathcal L}(F_r)\asymp r^{1-s}$ as $r \to 0$ then 
 $d(A, H, D) = \dimb F =s$ and
\be 
 \frac{ \sum_{k=1}^N \sigma_k (|D|^{-s})}{\log N} \asymp 1 \quad \mbox{ as } N\to \infty. \label{asyam}
\ee
In particular for any Dixmier trace
$ 0 < \tau_\omega (|D|^{-s})<\infty$ 
 and
\be
\tau_\omega (\xi(f)|D|^{-s}) = c\int_F f d\nu \label{ncint}
 \ee
defines a non-degenerate integral for some probability measure $\nu$ on $F$ where $0<c<\infty$.
Moreover, if $F$ is Minkowski measurable with Minkowski content $M_{s}(F)$ then $ \lvert D \rvert^{-s}$ is a measurable operator with
\be
-\hspace{-0.45cm}\int  \rvert D \lvert^{-s} =  \tau_\omega (|D|^{-s})= 2^{s}(1 - s)M_{s}(F). \label{cval}
\ee
\end{theo}

\begin{pf}
It is easy to verify that 
$$d(A, H, D)  = \inf\big\{s: \sum_{n=1}^\infty  | I_{n} |^s < \infty\big\} = \dimub F;$$
 the first equality follows from the definitions (see the proof of Corollary \ref{betasum} for a similar argument), and the second may be established using Minkowski's definition of upper box-counting dimension.   If $F$ is almost Minkowski measurable then $ | I_{n} |\asymp n^{-s}$, see \cite{GI1,LandP}, which gives (\ref{asyam}).
We conclude that each Dixmier trace lies between the lower and upper limits of (\ref{asyam}), and the trace state  $f \mapsto \tau_\omega (\xi(f)|D|^{-s})$ defines a Radon measure and integral on $F$. In the Minkowski measurable case the constant in (\ref{cval}) comes from the asymptotic limit $ | I_{n} |\sim 2^{s}(1 - s)M_{s}(F)n^{-s}$, see \cite{GI1,LandP} 
\end{pf}

Whilst self-similar fractals are not necessarily Minkowski measurable we nevertheless get a noncommutative integral.
Let  $\{ \psi_1,\ldots,\psi_m\}$ be an {\it iterated function system} where each $\psi_i$ is a contracting similarity on $\bbbr$ with ratio $r_i$, such that the intervals $\psi_i  [0,1]$ are disjoint and ordered, that is $0= \psi_1(0) <  \psi_1(1)<  \psi_2(0)< \ldots <  \psi_m(0)<  \psi_m(1)=1$. The {\it attractor} $F$ of the system is characterised as the unique non-empty compact set satisfying
\be F = \bigcup_{i=1}^{m} \psi_i(F), \label{ifs}
\ee
see \cite{F1}, such a set is called {\it self-similar}. Write $\epsilon_1, \ldots,\epsilon_{m-1}$ for the gap lengths between consecutive  intervals $\psi_i [0,1]$, so that $ \epsilon_i = \psi_{i+1}(0) -  \psi_{i}(1)$.

\begin{theo}\label{thm2.2}
Let $F$ be a self-similar subset of $[0, 1]$ with ratios $r_1,\ldots,r_m$ and gap lengths $ \epsilon_1,\ldots, \epsilon_{m-1}$. Then  $d(A, H, D)=\dimh F =\dimb F=s$  where $s$ is the positive solution of 
$\sum_{i=1}^m r_i^s=1$.  Moreover  the operator $\xi(f) \lvert D \rvert^{-s}$ is measurable for all $f \in C(F)$ with 
\be
\tau_\omega( \xi(f) \lvert D \rvert^{-s}) =-\hspace{-0.45cm}\int \xi(f) \rvert D \lvert^{-s} = c \int_F  f d\nu,
\label{tauss1} 
\ee
where $\nu$ is Hausdorff measure normalized so that $\nu(F)=1$,
and 
\be
c= \frac{2 \sum_{i = 1}^{m - 1}  \epsilon_i^{s}  }{
\sum_{i = 1}^{m } r_i^{s}\log (1/r_i)}
= \frac{2 \sum_{i = 1}^{m - 1}  \epsilon_i^{s}  }{ h_{\mu}},\label{tauss2} 
\ee
where $h_{\mu}$ is the measure theoretic entropy.
\end{theo}

\begin{pf}
Every self-similar fractal $F$  is almost Minkowski measurable with $d(A, H, D)=\dimh F =\dimb F$.  Thus $0<\tau_\omega( |D|^{-s})<\infty$ for any Dixmier trace by Theorem \ref{tracemm}.
A renewal theory argument (see Theorem \ref{thm3.5} below) shows that $|D|^{-s }$ is measurable and gives the stated value of $c= \tau_\omega( \lvert D \rvert^{-s})$. Then for each Dixmier trace $\tau_\omega( \xi(f) \lvert D \rvert^{-s})=c \int_F  f d\nu$ defines an integral for some probability measure $\nu$. A scaling argument shows that $\nu$ satisfies
$\nu(A) = \sum_{i=1}^m r_i^s  \nu (\psi^{-1}_i(A))$ which implies that $\nu$ is normalized $s$-dimensional Hausdorff measure, and also that $\xi({\bf 1}_I) \lvert D \rvert^{-s}$ is measurable for every interval $I= \psi_{i_1}\circ\psi_{i_2} \circ\cdots\circ \psi_{i_k}[0,1]$. It follows that $\xi(f) \lvert D \rvert^{-s}$ is measurable for all continuous $f$ and we can write the noncommutative integral in (\ref{tauss1}).
\end{pf}


\section{Spectral triples for multifractals}
\setcounter{equation}{0}
\setcounter{theo}{0}
In this section we seek analogues of the fractal results of Section 3 for multifractal measures. We present an approach for encoding the coarse multifractal spectrum of measures on an interval in terms of  spectral dimensions; in some cases, in particular for self-similar measures, the Dixmier trace corresponds to integration with respect to a multifractal auxillary measure.

Coarse multifractal analysis concerns the asymptotic behaviour of moment sums of a measure. Let $\mu$ be a Borel probability measure on $[0,1]$ which we will always assume has no atoms.
Let   $\mathcal{B}_{r}$  be denote the grid intervals of length $r$, that is
  $\mathcal{B}_{r}=\{[mr,(m+1)r] : m \in \mathbb{Z}\}$. 
We consider the behaviour of the moment sums  $\sum_{B \in \mathcal{B}_{r} } \mu(B)^{q}$ as $r \to 0$  and in particular limits such as
\be
\beta_0(q) =  \limsup_{r \to 0}\frac{ \log \sum_{B \in \mathcal{B}_{r} }\mu(B)^{q}}{-\log r} \label{defbeta0}
\ee 
The Legendre transform of $\beta_0$ is often known as the coarse multifractal spectrum of $\mu$ and under certain conditions equals the fine spectrum of $\mu$ which is defined in terms of the Hausdorff dimensions of the sets at which $\mu$ has given local dimension, see \cite{F1,MF}. 
Whilst the definition (\ref{defbeta0}) of $\beta_0(q) $ may be reasonable for $ q\geq 0$, the moment sums can become unstable for $q<0$ since for certain values of $r$ some of the intervals $B \in \mathcal{B}_{r}$ may just clip the support of the measure making $ \mu(B)^{q}$ uncharacteristically large. Various devices have been proposed to overcome this difficulty and our approach will also deal with negative $q$ in a stable manner.

We first recall a modification of (\ref{defbeta0}) introduced by Riedi \cite{Rie} that is more satisfactory for $q<0$. Given a  closed interval $B = [s,t] \subseteq\bbbr$ and $a\geq 0$, let ${\overline B}^a \equiv [s - \frac{1}{2}a(t-s), t+\frac{1}{2}a(t-s)] $ be the interval obtained by enlarging $B$ by a factor $a+1$ about its centre. We also let 
  $\mathcal{B}_{r}^* = \{B \in\mathcal{B}_{r}: \mu(B)>0\}$,   that is the intervals of  $\mathcal{B}_{r}$ with interiors that intersect the support of $\mu$.
  
Following Reidi \cite{Rie} we work with moment sums over these enlarged intervals and define 
\be
\beta(q) = \inf \Big\{\beta : \limsup_{r \to 0} r^\beta\sum_{B \in \mathcal{B}_{r}^* } \mu({\overline B}^a)^{q}\ = 0\Big\}
= \sup \Big\{\beta : \limsup_{r \to 0} r^\beta\sum_{B \in \mathcal{B}_{r}^* } \mu({\overline B}^a)^{q}\ = \infty\Big\}. \label{beta3}
\ee
(as is usual with such definitions the critical value $\beta(q)$ is where a limit of a sum jumps from $\infty$ to $0$).  
We see from the next proposition that this definition of $\beta(q)$ is independent of $a>0$ for all $q$, and independent of $a\geq 0$ for $q \geq 0$.
  
 \begin{prop}\label{newbeta}
Let $\mu$ be a probability measure on $[0,1]$. 

$(a)$ Given $q \geq 0$ and $0\leq a\leq b$ there is a constant $c_1>0$ such that for sufficiently small $r$

\be
 c_1\sum_{B \in \mathcal{B}_{r}^* }\mu({\BB}^b)^{q}
\leq \sum_{B \in \mathcal{B}_{r}^* }\mu({\BB}^a)^{q}
\leq \sum_{B \in \mathcal{B}_{r}^* }\mu({\BB}^b)^{q}.\label{propa}
\ee

$(b)$ Given  $q< 0$ and $0< a \leq b$ there is a constant $c_2>0$ such that for sufficiently small $r$
\be
\sum_{B \in \mathcal{B}_{r}^* }\mu({\BB}^b)^{q}
\leq \sum_{B \in \mathcal{B}_{r}^* }\mu({\BB}^a)^{q}
\leq c_2 \sum_{B \in \mathcal{B}_{ra/(2+b)}^* }\mu({\BB}^b)^{q}.\label{propb}
\ee
Thus the asymptotic behaviour of $\sum_{B \in \mathcal{B}_{r}^* }\mu({\BB}^a)^{q}$ as $r \to 0$ is independent of $a\geq 0$ if $q\geq 0$ and is essentially independent of $a>0$ for all $q$.
 In particular $\beta(q)$ given by  $(\ref{beta3})$ does not depend on the value of $a>0$ for all $q$,  and    also $\beta(q)=\beta_0(q)$  if $q\geq 0$.
 
\end{prop}

\begin{pf} This is essentially given in  \cite{Rie}. The inequalities (\ref{propa}) and (\ref{propb})are easily established using that $\BB^a \subseteq \BB^b$ if $a \leq b$ together with the fact that, given $r,r'>0$,
for each $B \in \mathcal{B}_{r}^*$ a  bounded number of intervals in $ \mathcal{B}_{r'}^* $ intersect $B$. The asymptotic properties are an immediate consequence. 
\end{pf}

We next assume that the support of the measure $\mu$ is a fractal $F$ of the type considered in Section 3.  Thus $F$ is a fractal subset of $[0, 1]$ of Lebesgue measure $0$ that is the complement in $[0, 1]$ of the intervals $ \{I_n\}$ ordered by decreasing length. The set  $F$ is said to satisfy the \textit{lacunarity condition}  with \textit{lacunarity constant}  $0<\lambda \leq 1$ if for all $x \in F$ and $0 < r \leq 1$ the interval $B(x, r)$ (centred at $x$ with length $2r$) contains a complementary interval $I_n$  of length greater then or equal to $\lambda r$  (such an $F$ is sometimes called \textit{porus}). Many sets, including all self-similar and self-conformal subsets of $[0,1]$ satisfy this condition.

The lacunarity condition imposes bounds on the dimensions of the sets and the rate of decrease of the complementary interval lengths.

\begin{lem}\label{lac}
Let $F \subseteq [0, 1]$  satisfy the lacunarity condition. Then 
\be
0 < \dimh F \leq \dimlb F \leq\dimub F \leq 1\label{lacdims}
\ee
and 
\be
 -\infty < \liminf_{n\to\infty} \frac{\log | I_{n} |}{\log n}
 \leq \limsup_{n\to\infty} \frac{\log | I_{n} |}{\log n} \leq -1.\label{laclengths}
 \ee
\end{lem}

\begin{pf}
Using the lacunarity condition repeatedly to replace single points by nearby pairs of points one may construct a subset of $F$ that is a Cantor set of positive Hausdorff dimension, so $0 < \dimh F$. The other inequalities of (\ref{lacdims}) are standard.

The bounds (\ref{laclengths}) are obtained in \cite[Proposition 3.7]{F2} (in fact with more precise bounds which depend on the lower and upper box dimensions of $F$).
\end{pf}

Keeping the above notation, note that $\IB^a = [b_{n}^{-}-\frac{1}{2}a(b_{n}^{+}-b_{n}^{-}), b_{n}^{+}+\frac{1}{2}a(b_{n}^{+}-b_{n}^{-})]$ is the enlargement of $I_n$ by a factor $a+1$ about its centre. Working with the enlarged intervals $\IB^a$ achieves two things: it enables us to relate  the moment sums to sums involving  the complementary intervals, and it gives stable behaviour when $q<0$.

\begin{prop}\label{propA}
Let $\mu$ be a probability measure with support the whole of $F$ where $F\subseteq [0,1]$ satisfies the lacunarity condition with constant $\lambda$. Let  $a \geq 2/\lambda$ and $b >0$. Then there are numbers $\eta_1, \eta_2,c_1,c_2>0$,  which depend on $a,b,\lambda$ and $q$ such that for sufficiently small $r$ 
\begin{equation}\label{gb}
c_1 \sum_{\lambda\eta_1 r \leq |I_{n}| \leq \eta_1 r} \mu(\IB^a)^{q}  \leq \sum_{B \in \mathcal{B}_{r}^* }\mu(\BB^b)^{q}
\leq c_2 \sum_{\lambda\eta_2 r \leq |I_{n}| \leq \eta_2 r} \mu(\IB^a)^{q}.
\end{equation}
\end{prop}

\begin{pf}
First fix $0<r< 1$ and write  ${\cal I} =\{I_n: \lambda r \leq |I_n| \leq r\}$. 
Let $x \in F$. By the lacunarity condition there is an interval $I_n=[b_{n}^{-}, b_{n}^{+}] \in {\cal I}$ such that 
$I_n \subseteq B(x,r)$. Thus 
\be
x \in  [b_{n}^{-}-r,  b_{n}^{+}+r]
\subseteq [b_{n}^{-}-{\textstyle \frac{1}{2}} a|I_n|,  b_{n}^{+}+{\textstyle \frac{1}{2}}a|I_n|] 
=  \IB^a   ,\label{inc1}
\ee
since $r \leq |I_n|/\lambda \leq \frac{1}{2} a|I_n|$.

On the other hand given $x\in[0,1]$, if $x \in   \IB^a=[b_{n}^{-}-{\textstyle \frac{1}{2}} a|I_n|,  b_{n}^{+}+{\textstyle \frac{1}{2}}a|I_n|] 
$ where $I_n \in {\cal I}$, then
$$ I_n\subseteq [x-({\textstyle \frac{1}{2}}a+1)|I_n|,x+({\textstyle \frac{1}{2}}a+1)|I_n|]
\subseteq [x-({\textstyle \frac{1}{2}}a+1)r,x+({\textstyle \frac{1}{2}}a+1)r].$$
Since $|I_n| \geq \lambda r$ and the intervals $|I_n|$ have disjoint interiors, $x$ lies in at most $(a+2)/\lambda$  intervals $\IB^a$ with $I_n \in {\cal I}$. Together with (\ref{inc1}) we see that each $x \in F$ lies in at least one and at most $(a+2)/\lambda$ of the intervals
 $\{\IB^a: I_n \in {\cal I}\}$.
 
The inequalities (\ref{gb}) may now be established using slightly different geometrical arguments for $q \geq0$ and $q<0$.

(a) For sufficiently small $r$ let $I_n$ be an interval such that  
$\lambda  br/((2+a)\leq |I_n| \leq br/((2+a)$.  Since the endpoints of $I_n$ are in the support of $\mu$ simple geometry shows that we can find 
$B \in \mathcal{B}_{r}^* $ such that $\IB^a\subseteq\BB^b $.
This gives the left hand inequality of (\ref{gb}) for  $q\geq 0$ with $\eta_1 =b/((2+a)$,  noting that the number of $ \BB^b$ with $B \in \mathcal{B}_{r}^* $ that contain such an $\IB^a$  is bounded independently of $r$. 

(b) For sufficiently small $r$ let $I_n$ be an interval such that  
$ r(2+b)/a \leq |I_n| \leq \lambda r(2+b)/a $.  Since the endpoints of $I_n$ are in the support of $\mu$, we see that there exists 
$B \in \mathcal{B}_{r}^* $ such that $\BB^b\subseteq \IB^a$.
This gives the left hand inequality of (\ref{gb}) for  $q< 0$ with $\eta_1 =\lambda(2+b)/a$, noting from above that each $\BB^b$ with $B \in \mathcal{B}_{r}^* $ is contained in a bounded number of different $\IB^a$ with $|I_n|$ in this range.
 
(c) Let $B \in \mathcal{B}_{r}^* $  for sufficiently small $r$.  Using the lacunarity condition, a geometrical estimate shows that there is an $I_n$ such that $\BB^b \subseteq\IB^a$ with 
 $\lambda(2+b)(2+a -2/\lambda)^{-1} r \leq |I_n| \leq (2+b)(2+a -2/\lambda)^{-1}r$. Taking powers of measures  and summing gives the right hand inequality of (\ref{gb}) for  $q\geq 0$ with $\eta_2 =(2+b)(2+a -2/\lambda)^{-1}$, noting from above that the number of $\IB^a$ with $|I_n|$ in this range that can contain each $ \BB^b$ is bounded independently of $r$.
 
(d) Let $B \in \mathcal{B}_{r}^* $  for sufficiently small $r$.  Using the lacunarity condition, we see  that there is an $I_n$ such that $\IB^a \subseteq \BB^b$ with 
 $\lambda b(2+a)^{-1} r \leq |I_n| \leq b(2+a)^{-1}r$. This gives the right hand inequality of (\ref{gb}) for  $q<0$ with $\eta_2 =b(2+a)^{-1}$, noting from above that the number of $\IB^a$ in this range contained in each $ \BB^b$ is bounded independently of $r$.
  
\end{pf}

The moment function $\beta(q)$ given by (\ref{beta3}) may now be expressed in terms of the extended complementary intervals.

\begin{prop}\label{betasum}
Let $\mu$ be a probability measure with support the whole of $F$ where $F\subseteq [0,1]$ satisfies the lacunarity condition. Let $a \geq 2/\lambda$. Then 
\begin{eqnarray}\label{beta1}
\beta(q)&=&
 \inf \left\{\beta : \sum_{n=1}^{\infty}  \mu(\IB^a)^{q}
|I_n|^\beta< \infty\right\} \\ 
&=&\inf\left\{\beta:
\limsup_{N \to \infty} \frac{ \sum_{n=1}^{N}  \mu(\IB^a)^{q}
|I_n|^\beta}{\log N}= 0\right\}.\label{beta2}
\end{eqnarray}
\end{prop}

\begin{pf}
Combining  (\ref{gb}) with definition  (\ref{beta3}) easily  gives (\ref{beta1}).

If  the series in (\ref{beta1}) converges then the numerator in  (\ref{beta2}) is bounded so the upper limit is $0$.
On the other hand, if $ \sum_{n=1}^{N}  \mu(\IB^a)^{q}|I_n|^\beta \leq c \log N$ for all $N$ then, 
for $\epsilon > 0$,$$ \sum_{n=2^k}^{2^{k+1}-1}  \mu(\IB^a)^{q}|I_n|^{\beta+\epsilon}
 \leq c (\log 2^{k+1})(2^{-k})^{\epsilon} \leq c_1 k 2^{-k \epsilon},$$
 since $|I_n| \leq 2^{-k}$ if $n\geq 2^k$, so summing over $k$ the series in  (\ref{beta1}) converges when  the exponent of $|I_n|$ equals $\beta + \epsilon$ for all $\epsilon >0$. 
\end{pf}

Note that  (\ref{beta1})  defines a $\beta(q)$ function even without the lacunarity condition or with $a \geq 2/\lambda$. However, the definition may then fail to detect the entire support of the measure. 

We now define functions $g$ and operators $G$  to express the multifractal behaviour in terms of a spectral triple and Dixmier trace.
Given $a> 0$ and a probability measure $\mu$ with support $F$ let
$g\equiv g_{a,\mu}  : \bigcup_{n =1}^\infty\{ b_{n}^{-}, b_{n}^{+}\} \rightarrow [0,\infty)$
be given by
\be
g(b_{n}^{-}) = \;g(b_{n}^{+})=\mu(\IB^a)>0. \label{gdef}
\ee
Thus $g$ associates with each interval $I_n$ the measure on either side of the interval that lies within a distance comparable to $|I_n|$.
To express this as an operator on $H$ we define $G\equiv G_{a,\mu}: H \rightarrow H$ to  be
\begin{equation}\label{Fmudelta}
G \mathrel{:=} \bigoplus_{n \in \mathbb{N}} \left( \begin{array}{cc} g(b_{n}^{-}) & 0 \\ 0 & g(b_{n}^{+}) \end{array} \right).
\end{equation}
(We avoid writing $ \xi(g)$ for $G$ since $g$ is not continuous.)

 The next proposition, which is a multifractal analogue of Theorem \ref{tracemm}, shows that the Dixmier traces
$\tau_\omega (G^q|D|^{-\beta}))$ reflect the behaviour of the multifractal moment sums when $F$ satisfies the lacunarity condition and $\mu$ satisfies a very mild density condition (\ref{mudensity}). Note that the condition (\ref{betaqexists}) is a multifractal analogue of a fractal being almost Minkowski measurable, a condition which is satisfied by many measures, including self-similar measures and Gibbs measures on cookie-cutter sets, with $\beta(q)$  the usual multifractal moment function given (for $q$ both positive and negative) by an indicial or pressure equation. Recall that $\{ \sigma_k (T) \}_{k=1}^{\infty}$ are the singular values of an operator $T$ arranged in decreasing order; thus
to find Dixmier  traces we have to reorder the terms of the sums by decreasing $g(b_n^{\pm})|I_n|$ rather than by decreasing  $|I_n|$.

\begin{theo}\label{multifracnoncomint}
Let $F \subseteq [0,1]$ satisfy the lacunarity condition with lacunarity constant $\lambda$ and suppose that the probability measure $\mu$ has support  $F$ and satisfies
\be
0< c_1 \leq\frac{\log\mu(B(x,r))}{\log r} \leq c_2 <\infty \quad (x \in F,\, 0<r\leq \rho) \label{mudensity}
\ee
for  some $c_1,c_2$ and $\rho>0$. Then, for $a \geq 2/\lambda$ and with $G\equiv G_{a,\mu}$,  
\be
\inf\{\beta:\tau_\omega (G^q|D|^{-\beta}) \in \mathcal{L}^{1,\infty}(H)\}=\beta(q), \label{betadixid}
\ee 
where $\beta(q)$ is given by $(\ref{beta3})$. 
If $\mu$ satisfies the condition
\be
\sum_{|I_{n}| \geq \tau}  \mu(\IB^a)^{q} \lvert I_{n} \rvert^{\beta} \asymp -\log \tau
\label{sumtau}
\ee
for small $\tau$, then $ \beta =\beta(q)$ and
\be
 \frac{ \sum_{k=1}^N \sigma_k (G^q|D|^{-\beta})}{\log N} \asymp 1 \quad \mbox{ as } N\to \infty. \label{posfihn}
\ee
In particular
$ 0 < \tau_\omega (G^q|D|^{-\beta}))<\infty$  for any Dixmier trace $\tau_\omega$, and
\be
  \tau_\omega (\xi(f)G^q|D|^{-\beta})) = \int_F f d\nu \label{integralf}
  \ee
defines a non-degenerate integral for some measure $\nu$.
Note  that $(\ref{sumtau})$ holds by Proposition \ref{propA} if  $\mu$ satisfies the moment condition
\be
\sum_{B \in \mathcal{B}_{r}^* }\mu(\BB^b)^{q} r^{\beta}\asymp 1 \quad (0<r \leq 1) \label{betaqexists}
\ee
for any $b>0$ (or for any $b\geq 0$ if $q\geq 0$).
\end{theo}

\begin{pf}
Note that the singular values $ \sigma_k (G^q|D|^{-\beta})$ are given by $\{\mu(\IB^a)^{q}| I_{n} |^{\beta}: n \in \mathbb{Z}^+\}$ each with multiplicity two, but to estimate the Dixmier trace we need to consider these in decreasing order.
For each complementary interval $I_n = [b_n^-, b_n^+]$ we have 
$$\big(B(b_{n}^{-},{\textstyle\frac{1}{2}}a| I_{n} |) \cup B(b_{n}^{+},{\textstyle\frac{1}{2}}a| I_{n} |) \big)\cap F=  \IB^a \cap F$$
so it follows from  (\ref{mudensity}) and (\ref{laclengths}) that
\be
0< c_3 \leq  \frac{\log   \mu(\IB^a)}{\log n} \leq c_4 <\infty \label{gdensity}
\ee
for large $n$, for suitable $c_3,c_4$.
Write
\be
{\cal N}(\tau) = \{n \in \mathbb{N} : \mu(\IB^a)^{q}|I_n|^\beta \geq \tau\}.\label{btau}
\ee
Using (\ref{laclengths}) and (\ref{gdensity})  there are numbers $0<m,p<\infty$ such that
\be
\{n: |I_n| \geq \tau^m\} \subseteq {\cal N}(\tau)  \subseteq\{n: |I_n| \geq \tau^p\} \label{nset}
\ee
for $\tau$ sufficiently small. Thus
$$\sum_{|I_{n}| \geq \tau^{m}}\mu(\IB^a)^{q}| I_{n} |^{\beta} 
\leq \sum_{  {\cal N}(\tau)}  \mu(\IB^a)^{q}| I_{n} |^{\beta} 
\leq \sum_{|I_{n}| \geq \tau^p}  \mu(\IB^a)^{q}| I_{n} |^{\beta}.
$$
From (\ref{laclengths}) $\log \mbox{card}(|I_{n}| \geq \tau) \asymp -\log\tau \asymp \log \mbox{card}(|I_{n}| \geq \tau^m) \asymp \log \mbox{card}(|I_{n}| \geq \tau^p)$, so for some $c_5,c_6 >0$,
\be
c_5 \frac{ \sum_{|I_{n}| \geq \tau^m} \mu(\IB^a)^{q}| I_{n} |^{\beta}}
{\log \mbox{card}(|I_{n}| \geq \tau^m)} \leq 
\frac{\sum_{ {\cal N}(\tau)}  \mu(\IB^a)^{q}| I_{n} |^{\beta} }{\log \mbox{card}({\cal N}(\tau))}
\leq c_6  \frac{ \sum_{|I_{n}| \geq \tau^p} \mu(\IB^a)^{q}| I_{n} |^{\beta}}
{\log \mbox{card}(|I_{n}| \geq \tau^p)}.\label{ratios}
\ee
for small $\tau$. Relating this to  (\ref{beta2})  and using (\ref{dixid}) gives (\ref{betadixid}).

Now suppose that (\ref{sumtau}) holds. By (\ref{laclengths}) $\log \mbox{card}(|I_{n}| \geq \tau^m) \asymp-\log\tau $,
so as $N \to \infty$, or equivalently as $\tau \to 0$, (\ref{nset}) and (\ref{ratios}) imply
$$ \frac{ \sum_{k=1}^N \sigma_k (G^q|D|^{-\beta})}{\log N} \asymp  \frac{\sum_{ {\cal N}(\tau)}  \mu(\IB^a)^{q}| I_{n} |^{\beta} }{\log \mbox{card}({\cal N}(\tau))} \asymp 1,$$
giving (\ref{posfihn}).

We conclude that each Dixmier trace lies between the implied lower and upper limits of (\ref{posfihn}), and the trace  $f \mapsto \tau_\omega (\xi(f)G^q|D|^{-\beta}))$ defines a Radon measure and integral on $F$.
\end{pf}

We now specialise to the case where $\mu$ is a self-similar measure on a self-similar subset $F$ of $[0,1]$. Let $\psi_i $ be contractions with ratios $r_i \,(i=1,\ldots,m) $ and let $F$ be as at  (\ref{ifs}); it is easy to see that $F$ satisfies the lacunarity condition. 
Let $E_1, \ldots,E_{m-1}$ be the gaps between the consecutive basic intervals $\psi_i [0,1]$, so that $E_i = [\psi_{i}(1),\psi_{i+1}(0)]$ and $|E_i| = \epsilon_i$; of course the $E_i$ will be amongst the $I_n$.
Given  $p_1,\ldots,p_m>0$ with $\sum_{i=1}^m p_i =1$  we take $\mu$ to be the unique probability measure satisfying 
\be
\mu(A) = \sum_{i=1}^m p_i \mu (\psi^{-1}_i(A)) \label{probrats}
\ee 
for all Borel sets $A$ so $\mu$ is supported by $F$; such a measure is a {\it self-similar measure}.

We define a number $r_0 \equiv r_0(a,\mu)$ by
\be
r_0 = \sum_{i=1}^{m}\sum_{I_n \subseteq \psi_i[0,1]}
|I_n|^{\beta}\big( \mu(\IB^a)^{q}-p_i^q  \mu(\psi_i^{-1}(\IB^a))^{q}\big)
+ \sum_{ i=1}^{m-1} |E_i|^{\beta} \mu({\overline E}_i^a)^{q}.
\label{defr}
\ee
We may think of $r_0$ as representing the difference between the terms is $(\ref{beta1})$ corresponding to $I_n$ and the scaled terms corresponding to  $ \psi _{i}^{-1}(I_n)$. Note in particular that the sums in (\ref{defr}) are finite sums: if $a|I_n|$  is less than the length of the smallest gap  $E_i$ the summands corresponding to $I_n$ in the left hand term will be $0$.

Here is the multifractal analogue of Theorem \ref{thm2.2}, which reduces to that theorem when $q=0$.

\begin{theo}\label{thm3.5}
Let $F$ be a self-similar subset of $[0, 1]$ as in $(\ref{ifs})$ and $\mu$ be the self-similar probability measure on $F$   given by $(\ref{probrats})$, with notation as above.  Then  $\beta(q)$  defined by $(\ref{beta3})$ is the solution of 
$\sum_{i=1}^m  p_i^q r_i^{\beta} =1$. Moreover the operator $ \xi(f) G^q \lvert D \rvert^{-\beta(q) }$ is measurable for all $f \in C(F)$ with\be
-\hspace{-0.45cm}\int \xi(f)G^q \lvert D \rvert^{-\beta(q) }=\tau_\omega( \xi(f) G^q \lvert D \rvert^{-\beta(q) }) =  c \int_F  f d\nu, \label{taumf1} 
\ee
where $\nu$ is the probability measure  on $F$  satisfying
\be
\nu(A) = \sum_{i=1}^m p_i^q r_i^\beta \nu (\psi^{-1}_i(A)) \label{nucond}
\ee 
and
\be
c = \frac{2 r_0}{\sum_{i = 1}^{m } p_i^q r_i^{\beta}\log (1/(p_i^q r_i^\beta))}\label{taumf2} 
\ee
where $r_0$ is given by $(\ref{defr})$.
\end{theo}

\begin{pf}
For such self-similar measures, $\beta(q)$ is the solution of  $\sum_{i=1}^m  p_i^q r_i^{\beta} =1$ and the moment condition (\ref{betaqexists}) satisfied, see \cite{F1,Rie}. Thus $0<\tau_\omega( G^q \lvert D \rvert^{-\beta(q) })<\infty$ for every Dixmier trace, by Theorem \ref{multifracnoncomint}.

To show that $G^q \lvert D \rvert^{-\beta(q) }$ is measurable and to find the constant we use a renewal theory argument. 
Let $\sigma_n = \mu(\IB^a)^{q} |I_n|^\beta = g(b_n^-)^q |I_n|^\beta  = g(b_n^-)^q |I_n|^\beta$ be the pair of  singular values of $G^q \lvert D \rvert^{-\beta(q)}$ associated with the interval $I_n$. 
For $\tau >0$ let 
$$S(\tau)= \sum_{\sigma_n \geq \tau}\sigma_n
= \sum_{\sigma_n \geq \tau}\mu(\IB^a)^{q} |I_n|^\beta$$ 
and
$$
 r(\tau)
 = \sum_{i=1}^{m}\sum_{I_n \subseteq \psi_i[0,1], |I_n| \geq \tau}
|I_n|^{\beta}\big( \mu(\IB^a)^{q}-p_i^q  \mu(\psi_i^{-1}(\IB^a))^{q}\big)
+ \sum_{ |E_i| \geq \tau} |E_i|^{\beta} \mu({\overline E}_i^a)^{q}.
$$
Note that $r(\tau)$ is decreasing with $r(\tau) = r_0$ (see (\ref{defr})) if $\tau$ is sufficiently small.
By comparing the measures of the intervals $\IB^a$ and $\psi _{i}^{-1}\IB^a$,  we get the scaling relation
$$ S(\tau) =  \sum_{i=1}^{m}p_i^q r_i^\beta S(\tau/(p_i^q r_i^\beta)) + r(\tau).$$
Transforming by setting $t= -\log \tau$ and $h(t) = S(e^{-t})$ gives
$$ h(t) =  \sum_{i=1}^{m}p_i^q r_i^\beta h(t + \log(p_i^q r_i^\beta)) + r(e^{-t}).$$ 
With $\beta = \beta(q)$ we have $\sum_{i=1}^m  p_i^q r_i^{\beta} =1$ so this becomes a renewal equation and we may apply the (elementary) renewal theorem (see \cite{Fel}) to conclude that
\be
\lim_{\tau\to 0}\frac{S(\tau)}{-\log \tau}  
= \lim_{t\to \infty}\frac{h(t)}{t}  
= \frac{\lim_{t\to \infty}r(e^{-t})}{\sum_{i = 1}^{m }p_i^q r_i^{\beta}\log (1/(p_i^q r_i^\beta))}
=\frac{r_0}{\sum_{i = 1}^{m } p_i^q r_i^{\beta}\log (1/(p_i^q r_i^\beta))}, \label{lim1}
\ee
noting that $r(\tau) =0$ for $\tau \geq 1$ and $r(\tau) =  r_0$ for $\tau$ sufficiently small. 

Now let $S_1(\tau) =   \#\{n:  \sigma_n \geq \tau\} = \#\{n:  \mu(\IB^a)^{q} |I_n|^\beta\geq \tau\} $  and 
\begin{align*}
r_1(\tau) =\sum_{i=1}^m \#  \big\{n: I_n \subseteq \psi_i[0,1]  &  \mbox{ with }  |I_n|^\beta\mu(\IB^a)^{q}  \geq \tau \mbox{ but }  
 |I_n|^\beta p_i^q\mu(\psi _{i}^{-1}(\IB^a))^q < \tau\big\}\\
&+  \#\{i:  |E_i|^{\beta} \mu({\overline E}_i^a)^{q} \geq \tau\}.
 \end{align*}
Note that $r_1$ is decreasing, with $r_1(\tau)=0$ if $\tau \geq 1$ and $r_1(\tau)= r_1>0$, say, for $ \tau $ sufficiently small.
Then we get another scaling relationship
$$ S_1(\tau) =  \sum_{i=1}^{m} S_1(\tau/(p_i^q r_i^\beta)) + r_1(\tau).$$
Setting $t= -\log \tau$ and $h(t) = e^{-t}S_1(e^{-t})$ gives the renewal equation
$$ h(t) =  \sum_{i=1}^{n} p_i^q r_i^\beta h(t + \log p_i^q r_i^\beta) +e^{-t} r_1(e^{-t}).$$ 
The renewal theorem (see \cite{F2,Fel}) implies that
$h(t) \asymp 1$ as $t \to \infty$ so  $ S_1(\tau) \asymp \tau^{-1}$ as $\tau \to \infty$ (in fact $h(t)$ either converges to a constant or is asymptotic to a periodic function, but we do not need that here.) We conclude that  $\log \#\{n:  \sigma_n \geq \tau\} /- \log \tau \to 1$ as $\tau \to 0$. Combining this with (\ref{lim1}) and recalling the definition of the Dixmier trace (\ref{Dixtraa})
$$\tau_{\omega}(G^q \lvert D \rvert^{-\beta(q) })=
\lim_{\tau \to 0} \frac{2S(\tau)}{ \log (\#\{n:  \sigma_n \geq \tau\})}
=\frac{2 r_0}{\sum_{i = 1}^{m } p_i^q r_i^{\beta}\log (1/(p_i^q r_i^\beta))}=c,$$
so that  $G^q \lvert D \rvert^{-\beta(q) }$ is measurable and $-\hspace{-0.45cm}\int G^q \lvert D \rvert^{-\beta(q) }=c$. 

For each Dixmier trace $\tau_\omega( \xi(f) \lvert D \rvert^{-s})=c \int_F  f d\nu$ defines an integral for some probability measure $\nu$. A scaling argument now allows us to identify the measure $\nu$. Given a word $(i_1,i_2,\ldots,i_k)$ with $i_j \in\{1,2,\ldots,m\}$,
let $I$ be the interval $\psi_{i_1}\circ\psi_{i_2} \circ\cdots\circ \psi_{i_k} [0,1]$. Writing ${\bf 1}_I$ for the indicator function of $I$ we note that the singular values of  $\xi({\bf 1}_I) G^q \lvert D \rvert^{-\beta(q) }$ are just those corresponding to the complementary intervals $\{I_n\}$ that are contained in $I$, apart from endpoints of complementary intervals that abut $I$. With this finite collection of endpoints as exceptions (and a finite set of singular values does not affect the singular trace), the mapping $\psi_{\bf i} \equiv\psi_{i_1}\circ\psi_{i_2} \circ\cdots\circ \psi_{i_k}$ gives a bijection between $\{I_n: n \in \mathbb{Z}^+\}$ and $\{I_n: I_n \subseteq I\}$. Moreover, $|\psi_{\bf i} (J)| = r_{i_1}r_{i_2} \ldots r_{i_k}|J|$ for any interval $J$  and 
$\nu(\psi_{\bf i} (J)) = p_{i_1}p_{i_2} \ldots p_{i_k}\nu(J)$ for any sufficiently small interval $J$. Thus 
$$\tau_\omega ( \xi({\bf 1}_I) G^q \lvert D \rvert^{-\beta(q) }) =  (p_{i_1}p_{i_2} \ldots p_{i_k})^q
(r_{i_1}r_{i_2} \ldots r_{i_k})^{-\beta(q)} \tau_\omega ( G^q \lvert D \rvert^{-\beta(q) })
= \nu(I) c $$
where $\nu$ is the measure defined by (\ref{nucond}), so that $\xi({\bf 1}_I) G^q \lvert D \rvert^{-\beta(q) }$ is measurable.
By extension $\tau_\omega( \xi(f) G^q \lvert D \rvert^{-\beta(q) }) = c \int f d\mu$ with $\xi(f) G^q \lvert D \rvert^{-\beta(q) }$ measurable, giving the noncommutative integral (\ref{taumf1}).
\end{pf}


\section{Further remarks}
\setcounter{equation}{0}
\setcounter{theo}{0}

This approach to representing coarse multifractal properties  in terms of spectral triples and Dixmier traces could certainly be extended to other settings. For example, one could work with limit fractals \cite{GI1,GI2}, or  an analogue of Theorem 3.6 might be obtained for Gibbs measures on cookie-cutter sets  using non-linear versions of the renewal theorem, see \cite{Lally}. 

In our analysis $g(b_{n}^{-})=g(b_{n}^{+})$ encoded the measures of intervals that are enlargements of $I_n$ about its centre. In a similar, but notationally more complicated way, we could take 
$g(b_{n}^{-}) = \mu [b_{n}^{-}-\frac{1}{2}a(b_{n}^{+}-b_{n}^{-}), b_{n}^{-}]$ and
$g(b_{n}^{+}) = \mu [b_{n}^{+}, b_{n}^{+}+\frac{1}{2}a(b_{n}^{+}-b_{n}^{-})]$ so that  $b_{n}^{-}$ and  $b_{n}^{-}$ reflect the measure lying on each side of the $I_n$. We might then set one or other of $g(b_{n}^{-}) $ or $g(b_{n}^{+}) $ to $0$ to permit a  `one-sided' multifractal analysis, see \cite{F4}.

Various other  constructions of spectral triples on fractals have been proposed. For example, \cite{GI2} uses a sequence of discrete point sets that approximate the fractal,  \cite{STCS} uses covering intervals and  \cite{CIL} employs a construction built on curves. Such constructions extend to measures on higher dimensional spaces and multifractal versions might also be possible.

\end{document}